\documentclass[a4paper,11pt]{article}
\usepackage{amsmath,amsthm,amssymb}
\usepackage[mathscr]{eucal}

\numberwithin{equation}{section}

{\theoremstyle{definition}\newtheorem{definition}{Definition}[section]

\newtheorem{remark}[definition]{Remark}
}

\newtheorem{proposition}[definition]{Proposition}
\newtheorem{lemma}[definition]{Lemma}
\newtheorem{theorem}[definition]{Theorem}

\newcommand{\F}{\mathbb{F}}
\newcommand{\cR}{\mathcal{R}}
\newcommand{\actson}[1][]{\overset{#1}{\curvearrowright}}
\newcommand{\SL}{\operatorname{SL}}
\newcommand{\rL}{\operatorname{L}}
\newcommand{\Aut}{\operatorname{Aut}}
\newcommand{\Out}{\operatorname{Out}}
\newcommand{\N}{\mathbb{N}}
\newcommand{\T}{\mathbb{T}}
\newcommand{\Z}{\mathbb{Z}}
\newcommand{\cF}{\mathcal{F}}
\newcommand{\cA}{\mathcal{A}}
\newcommand{\cV}{\mathcal{V}}
\newcommand{\id}{\mathord{\operatorname{id}}}
\newcommand{\si}{\sigma}
\newcommand{\cE}{\mathcal{E}}
\newcommand{\recht}{\rightarrow}
\newcommand{\cU}{\mathcal{U}}
\newcommand{\vphi}{\varphi}
\newcommand{\cW}{\mathcal{W}}
\newcommand{\R}{\mathbb{R}}
\newcommand{\al}{\alpha}
\newcommand{\eps}{\varepsilon}
\newcommand{\module}{\operatorname{mod}}
\newcommand{\Tr}{\operatorname{Tr}}
\newcommand{\ovt}{\overline{\otimes}}
\newcommand{\B}{\operatorname{B}}
\newcommand{\paut}{\phi}
\newcommand{\bH}{\operatorname{H}}
\newcommand{\bB}{\operatorname{B}}
\newcommand{\bZ}{\operatorname{Z}}
\newcommand{\bG}{\mathbb{G}}
\newcommand{\om}{\omega}
\newcommand{\cP}{\mathcal{P}}

\newcommand{\Qp}{\mathbb{Q}^*_+}
\newcommand{\No}{\mathbb{N} \setminus \{0\}}
\newcommand{\Ker}{\operatorname{Ker}}
\newcommand{\cK}{\mathcal{K}}
\newcommand{\convex}{\operatorname{conv}}
\newcommand{\Gammatil}{\widetilde{\Gamma}}
\newcommand{\Emb}[1][F]{\operatorname{Emb}(H_{#1}^\infty,(H_{#1} * K)^\infty)}
\newcommand{\Embed}{\operatorname{Emb}}
\newcommand{\cH}{\mathcal{H}}
\newcommand{\Rp}{\R_+}
\newcommand{\cJ}{\mathcal{J}}
\newcommand{\ot}{\otimes}
\newcommand{\cL}{\mathcal{L}}
\newcommand{\invlimit}{\mathop{\underleftarrow{\operatorname{lim}}}}

\newcommand{\Lambdatil}{\widetilde{\Lambda}}

\newcommand{\Ad}{\operatorname{Ad}}
\newcommand{\Centr}{\operatorname{Centr}}
\newcommand{\cG}{\mathcal{G}}
\newcommand{\pZ}[1][p]{\frac{\Z}{#1 \Z}}

\newcommand{\cM}{\mathcal{M}}
\newcommand{\dpr}{^{\prime\prime}}

\newcommand{\bS}{\mathscr{S}}

\newcommand{\Sfactor}{\bS_{\text{\rm factor}}}
\newcommand{\Sequiv}{\bS_{\text{\rm equiv}}}
\newcommand{\Scentr}{\bS_{\text{\rm centr}}}

\begin{document}
\begin{center}
{\LARGE\bf Actions of $\F_\infty$ whose II$_1$ factors and orbit
equivalence \vspace{0.5ex}\\ relations have prescribed fundamental
group}

\bigskip

{\sc by Sorin Popa\footnote{Partially supported by NSF Grant
DMS-0601082}\footnote{Mathematics Department; University of
    California at Los Angeles, CA 90095-1555 (United States).
    \\ E-mail: popa@math.ucla.edu} and Stefaan Vaes\footnote{Partially
    supported by Research Programme G.0231.07 of the Research Foundation --
    Flanders (FWO) and the Marie Curie Research Training Network
    Non-Commutative Geometry MRTN-CT-2006-031962.}\footnote{Department of Mathematics;
    K.U.Leuven; Celestijnenlaan 200B; B--3001 Leuven (Belgium).
    \\ E-mail: stefaan.vaes@wis.kuleuven.be}}
\end{center}

\begin{abstract}
\noindent We show that given any subgroup $\cF$ of $\R_+$ which
is either countable or belongs to a certain ``large'' class
of uncountable subgroups, there exist continuously many
free ergodic measure preserving actions $\sigma_i$
of the free group with infinitely many generators $\F_\infty$
on probability measure spaces $(X_i,\mu_i)$ such that their associated
group measure space II$_1$ factors $M_i=\rL^\infty(X_i)\rtimes_{\sigma_i} \F_\infty$
and orbit equivalence relations $\cR_i=\cR(\F_\infty\actson X_i)$
have fundamental group equal to $\cF$ and with $M_i$ (respectively
$\cR_i$) stably non-isomorphic. Moreover, these actions can be taken
so that $\cR_i$ has no outer automorphisms and any automorphism of
$M_i$ is unitary conjugate to an automorphism that acts trivially on
$\rL^\infty(X_i)\subset M_i$.
\end{abstract}

\section{Introduction}

There has recently been increasing interest in the study of rigidity
properties of II$_1$ factors $M= \rL^\infty(X)\rtimes \Gamma$ and
orbit equivalence relations $\cR=\cR(\Gamma \actson X)$, arising
from ergodic measure preserving (m.p.) actions $\Gamma \actson X$ of
countable non-amenable groups $\Gamma$ on probability measure spaces
$(X,\mu)$, via the {\it group measure space construction} of Murray
and von Neumann (\cite{MvN1}) and its generalization in \cite{fm}.
Some of the most intriguing phenomena concern the {\it fundamental
group} $\cF(M)$, $\cF(\cR)$ of such II$_1$ factors and equivalence
relations (\cite{MvN2}). While a lot of progress has been made in
this direction, several basic questions on how these invariants
depend on the isomorphism class of the group $\Gamma$ and on the
nature of the action $\Gamma \actson X$ remained open. Many of them
are variations of Murray-von Neumann's long standing question, on
what subgroups of $\Rp$ may occur as fundamental groups of II$_1$
factors (\cite{MvN2}).

For instance, while any countable subgroup $\cF\subset \Rp$ was
shown to occur as the fundamental group of II$_1$ factors and
equivalence relations arising from  {\it non-free} ergodic m.p.\
actions in \cite{P-strong}, the problem of whether uncountable
groups $\neq \Rp$ can appear, remained wide open. If in addition the
action $\Gamma \actson X$ is required to be {\it free}, then much
less is known: the only groups shown to appear as $\cF(M), \cF(\cR)$
were $\Rp$ itself (\cite{MvN2}) and the trivial group $\{1\}$
(\cite{gg}, \cite{gab1}, \cite{gab2}, \cite{fur1}, for equivalence
relations; respectively \cite{P-Betti}, \cite{P-strong},
\cite{P-gap}, \cite{IPP} for II$_1$ factors). It has been speculated
that $\Gamma \actson X$ free may in fact force $\cF(M), \cF(\cR)$ to
be rational, when $\neq \Rp$. The case $\Gamma=\F_\infty$ is of
particular interest, due to its universality properties and to the
fact that, while $\cF(\cR)=\{1\}$ for any free ergodic m.p.\ action
$\F_n \actson X$, $2\leq n < \infty$ (by \cite{gab1}, \cite{gab2})
and some classes of
$\F_n$-actions are known to produce factors $M$ with
$\cF(M)=\{1\}$ (\cite{P-Betti}, \cite{OPI}), there was not one single case of a
free ergodic $\F_\infty$-action for which $\cF(\cR), \cF(M)$ could
be calculated.

We solve these problems here, by exhibiting a large family $\bS$ of
subgroups of $\R$, containing $\R$ itself, all of its countable
subgroups, as well as uncountable subgroups with any Hausdorff
dimension $\alpha \in (0,1)$, such that for each $H \in \bS$ there
exist ``many'' free ergodic m.p.\ actions of $\F_\infty$ for which
both the associated II$_1$ factor $M$ and orbit equivalence relation
$\cR$ have fundamental group equal to $\exp(H)$. Moreover, we
construct these actions so that we also have a good control on the
{\it outer automorphism group} of $M$ and $\cR$. In particular, we
obtain the first examples of II$_1$ factors with separable predual
(resp.\ standard countable m.p.\ ergodic equivalence relations) having
uncountable fundamental group different from $\Rp$.

The description of the set $\bS$ depends on results and notations
from \cite{aar-nad}, \cite{MNP}, which we recall for convenience.
Thus, following \cite{aar-nad} and \cite{MNP}, we call
\emph{ergodic measure on $\R$} any $\sigma$-finite measure $\nu$ on
the Borel sets of $\R$ satisfying the following properties, where we
denote $\lambda_x(y) = x+y$.
\begin{itemize}
\item For all $x \in \R$, either $\nu \circ \lambda_x = \nu$ or
$\nu \circ \lambda_x \perp \nu$.
\item There exists a countable subgroup $Q \subset \R$
such that $\nu \circ \lambda_x = \nu$ for all $x \in Q$
and such that every $Q$-invariant Borel function on $\R$ is
$\nu$-almost everywhere constant.
\end{itemize}
For every ergodic measure $\nu$ on $\R$, one defines
$$H_\nu := \{x \in \R \mid \nu \circ \lambda_x = \nu \} \; .$$
The obvious ergodic measures $\nu$ are the Lebesgue measure,
with $H_\nu = \R$, and the counting measure on a countable
subgroup $H \subset \R$, having $H_\nu = H$.
But, in \cite[Section 4]{aar-nad}, ergodic measures $\nu$
are constructed such that $H_\nu$ is an uncountable subgroup of $\R$
with prescribed Hausdorff dimension $\alpha \in (0,1)$. We'll recall this
construction in the preliminary section of the paper.

The family $\bS$ of subgroups $H$ of $\R$ for which we can
construct free ergodic m.p.\ actions of $\F_\infty$ whose II$_1$
factors and equivalence relations have fundamental group $\cF=\exp(H)$ will
consist of all subgroups of $\R$ of the form $H_\nu$. Thus,
we prove:

\begin{theorem} \label{main-theorem}
Let $\nu$ be an ergodic measure on $\R$ and $\cG$ any totally
disconnected unimodular locally compact group. There exists an
uncountable family $(\si_i)_{i \in I}$ of free ergodic m.p.\ actions
of the free group $\F_\infty$ on probability measure spaces,
$\si_i : \F_\infty \actson (X_i,\mu_i)$, such that the associated
orbit equivalence relations $\cR_i = \cR(\F_\infty \actson[\si_i]
X_i)$ and II$_1$ factors $M_i = \rL^\infty(X_i) \rtimes_{\si_i}
\F_\infty$ have the following properties.
\begin{itemize}
\item The fundamental groups are given by $\cF(\cR_i) = \cF(M_i) = \exp(H_\nu)$.
\item The outer automorphism groups are given by $\Out(\cR_i) \cong \cG$ and
$\Out(M_i) \cong \bH^1(\si_i) \rtimes \cG$. Here $\bH^1(\si_i)$ denotes
the $1$-cohomology group of $\si_i$
with values in $\T$ and corresponds to automorphisms of $M_i$ that are the
identity on $\rL^\infty(X_i)$.
\item The II$_1$ factors $M_i$ are not stably isomorphic and, in particular,
the equivalence relations $\cR_i$ are not stably orbit equivalent.
\end{itemize}
\end{theorem}

In particular, the above result provides a large class of free
ergodic m.p.\ actions of $\F_\infty$ whose orbit equivalence
relations $\cR$ have trivial outer automorphism group,
$\Out(\cR)=\{1\}$. Previous examples of free ergodic m.p.\ group-actions
$\Gamma\curvearrowright X$ with $\Out(\cR(\Gamma \actson X)) = \{1\}$ were
constructed in \cite{ge1},
\cite{fur2}, \cite{monod-shalom}, \cite{IPP}, \cite{PV},
\cite{vaes}, but for groups $\Gamma$ which either are higher rank
lattices (\cite{ge1}, \cite{fur2}), contain infinite subgroups with
the relative property (T) (\cite{IPP}, \cite{PV}, \cite{vaes}) or
are products of word-hyperbolic groups (\cite{monod-shalom}). It was
an open problem on whether actions of free groups can have this
property.

There is an interesting parallel between the above result, on groups
that may occur as the ``symmetry groups'' $\cF(M)$, $\Out(M)$ of
II$_1$ factors $M=\rL^\infty(X)\rtimes \F_\infty$ arising from free
ergodic m.p.\ actions of $\F_\infty$, and the results in \cite{PS},
showing that any group-like object that occurs in Jones' theory of
subfactors can be realized as ``generalized symmetries'' of
$M=L(\F_\infty)$. Along these lines, we conjecture that any group
$\cF\subset \R_+$ that can be realized as the fundamental group of a
separable II$_1$ factor $M$ (respectively as fundamental group of a
countable, m.p.\ ergodic equivalence relation $\cR$) can in fact be
realized as the fundamental group of a II$_1$ factor of the form
$M=\rL^\infty(X)\rtimes \F_\infty$ (resp.\ equivalence relation
$\cR(\F_\infty \actson X)$), for some appropriate free ergodic m.p.\
action $\F_\infty \actson X$.

The part concerning orbit equivalence relations of the Theorem \ref{main-theorem}
above is proved in Section~5, see Theorem \ref{thm.prescribed} and see also
Remark \ref{rem.extras} for a slight generalization.
Since the actions that we construct will be HT, in the sense of
\cite{P-Betti}, the part of Theorem \ref{main-theorem} concerning
the associated II$_1$ factors will follow automatically, by results
in \cite{P-Betti} (see Section~\ref{sec.vnalg}).

\subsubsection*{Acknowledgment}

The second author would like to thank the Mathematics Department
at UCLA for their warm hospitality during the work on this paper.

\section{Notations and preliminaries}

The fundamental group $\cF(M)$ of a II$_1$ factor $M$ is the
subgroup of $\Rp$ defined by
$$\cF(M) = \{ \tau(p)/\tau(q) \mid pMp \cong qMq \}$$
where $p,q$ are non-zero orthogonal projections in $M$. Similarly,
the fundamental group $\cF(\cR)$ of a countable measure preserving
(m.p.) ergodic equivalence relation $\cR$ on the standard
probability measure space $(X,\mu)$ (as defined in \cite{fm},
hereafter called {\it II$_1$  equivalence relation}) is given by
$$\cF(\cR) = \{ \mu(Y) / \mu(Z) \mid \cR|_Y \cong \cR|_Z \} \; .$$

One should point out that, while the fundamental group $\cF(M)$ is
defined for any II$_1$ factor $M$, we consider in this paper only
{\it separable} von Neumann factors, i.e. factors with separable
predual (equivalently, factors acting on separable Hilbert spaces).
For a II$_1$ factor, this is the same as requiring $M$ to be
separable in the Hilbert norm $\|x\|_2=\tau(x^*x)^{1/2}$, given by
the (unique) normalized trace $\tau$ on $M$, and is also equivalent
to $M$ being countably generated (as a von Neumann algebra). Thus,
any II$_1$ factor $M$ associated to a II$_1$ equivalence relation
$\cR$, via the generalized group measure space construction of
Feldman-Moore \cite{fm}, is separable, and one clearly has $\cF(\cR)
\subset \cF(M)$. This inclusion may be strict, in fact there even
exist free ergodic actions for which $\cF(\cR)=\{1\}$ while
$\cF(M)=\R_+$ (cf.\ 6.1 in \cite{P-gap}, based on \cite{CJ}).

\vskip .05in

Suppose now that $\Lambda \actson (X,\mu)$ is a free ergodic m.p.\
action. Define $(Y,\eta) = (X \times \Z, \mu \times c)$, where $c$
denotes the counting measure. We set $\Lambda^\infty := \Lambda
\times \Z$ acting freely on $(Y,\eta)$ as the product of the action
$\Lambda \actson X$ and the translation action $\Z \actson \Z$.
Recall in this respect that a non-singular automorphism $\Delta$ of
a measure space $(Y,\eta)$ is called (essentially) {\it free} if the
set $\{y \in Y \mid \Delta(y) = y \}$ has measure zero. Also, an
action of a countable group $\Gamma$ on $(Y,\eta)$, by non-singular
automorphisms, is free if each non-trivial element of the group
implements a free automorphism.
\begin{itemize}
\item We denote by $\cR(\Lambda \actson X)$ the orbit equivalence relation on $X$,
i.e.\ $x \sim y$ if and only if there exists $g \in \Gamma$ with $y = g \cdot x$.
We also consider the amplified equivalence relation $\cR(\Lambda^\infty \actson Y)$
and note that $(x,n) \sim (y,m)$ if and only if $x = g \cdot y$ for some $g \in \Lambda$,
\item We denote by $[\Lambda]$, resp.\ $[\Lambda^\infty]$ the group of non-singular
automorphisms $\Delta$ of $(X,\mu)$, resp.\ $(Y,\eta)$ satisfying $\Delta(x) \sim x$
almost everywhere. Note that the automorphisms in $[\Lambda]$ and $[\Lambda^\infty]$
are measure preserving.
\item We denote by $\Aut(\cR(\Lambda \actson X))$ the group of non-singular
automorphisms of $(X,\mu)$ satisfying $\Delta(x) \sim \Delta(y)$ if and only
$x \sim y$. We identify automorphisms that are equal almost everywhere.
We similarly consider $\Aut(\cR(\Lambda^\infty \actson Y))$.
An automorphism $\Delta \in \Aut(\cR(\Lambda \actson X))$
is automatically measure preserving, while an automorphism
$\Delta \in \Aut(\cR(\Lambda^\infty \actson Y))$ automatically
scales the measure $\eta$ with a positive constant that we denote
by $\module \Delta$.
\item Note that $[\Lambda], [\Lambda^\infty]$ are normal subgroups
of $\Aut(\cR(\Lambda \actson X))$, $\Aut(\cR(\Lambda^\infty \actson Y))$.
We denote by $\Out(\cR(\Lambda \actson X))$ and $\Out(\cR(\Lambda^\infty \actson Y))$
the corresponding quotient groups.
\item We identify a measure preserving automorphism $\Delta$ of $X$ with the
automorphism $\Delta \times \id$ of $Y$.
\item The maps $\Delta \mapsto \Delta \times \id$ and $\Delta \mapsto
\module \Delta$ induce an exact sequence
\begin{equation}\label{eq.exactsequence}
e \longrightarrow \Out(\cR(\Lambda \actson X)) \longrightarrow
\Out(\cR(\Lambda^\infty \actson Y)) \longrightarrow \cF(\cR(\Lambda \actson X))
\longrightarrow e \; .
\end{equation}
In particular, triviality of $\Out(\cR(\Lambda^\infty \actson Y))$ is
equivalent with triviality of both the outer automorphism group
$\Out(\cR(\Lambda \actson X))$ and the fundamental group
$\cF(\cR(\Lambda \actson X))$.
\item Suppose that $\Gamma < \Lambda$ is a subgroup and $\Gamma \actson X$
is still ergodic. We define $\Embed(\Gamma^\infty,\Lambda^\infty)$
as the set of non-singular automorphisms $\Delta$ of $Y$ such that
for all $g \in \Gamma^\infty$ and almost all $y \in Y$, we have
$\Delta(g \cdot y) \in \Lambda^\infty \cdot \Delta(y)$. Clearly,
$\Embed(\Gamma^\infty,\Lambda^\infty)$ is stable under composition
on the left by elements of $[\Lambda^\infty]$. We say that $\Delta$
is outer if $\Delta \not\in [\Lambda^\infty]$.
\end{itemize}

In the formulation of Theorem \ref{main-theorem} and in the last
section of this article, we use the $1$-cohomology group
$\bH^1(\si)$ of a free ergodic m.p.\  action $\si : \Lambda \actson
X$. Recall that $\bZ^1(\si)$ is the abelian Polish group of
functions $\om : \Lambda \times X \recht \T$ satisfying $\om(gh,x) =
\om(g,\si_h(x)) \om(h,x)$ almost everywhere. The subgroup
$\bB^1(\si)$ consists of $\om$ having the form $\om(g,x) =
\vphi(\si_g(x))\overline{\vphi(x)}$. The $1$-cohomology group
$\bH^1(\si)$ is by definition the quotient of $\bZ^1(\si)$ by
$\bB^1(\si)$ (\cite{singer}, \cite{fm}). This $\bH^1(\si)$ is a
stable orbit equivalence invariant (\cite{P-comp}, \cite{fm}).

\subsection*{Some families of subgroups of $\R$}

Since Theorem \ref{main-theorem} provides the first examples of II$_1$
factors and equivalence relations with uncountable fundamental group
different from $\R_+$, we are interested in the following sets of subgroups of the real line.

\begin{align*}
\Sfactor := \{ H \subset \R \mid \; & \text{there exists a separable
II$_1$ factor $M$ with}\; \cF(M) = \exp(H) \} \; , \\ \Sequiv := \{ H
\subset \R \mid \; & \text{there exists a II$_1$ equivalence
relation $\cR$ with}\; \cF(\cR) = \exp(H) \} \; .
\end{align*}
A consequence of our main result says that both $\Sfactor$ and $\Sequiv$ contain
$$\bS := \{ H_\nu \mid \nu \;\;\text{is an ergodic measure on $\R$} \} \; .$$
In fact, in the course of the proof of Theorem \ref{thm.prescribed}, we will see that
\begin{equation}\label{eq.maininclusion}
\bS \subset \Scentr \subset \; \Sfactor \cap \Sequiv \; ,
\end{equation}
where the set $\Scentr$ is defined as follows. Suppose that $(Y,\eta)$
is a standard infinite measure space and
$\Lambda \actson (Y,\eta)$ an ergodic, measure preserving action.
Define $\Centr_\Lambda(Y)$ as the subgroup of $\Aut(Y,\eta)$ consisting
of automorphisms $\Delta$ that commute with the $\Lambda$-action.
If $\Delta \in \Centr_\Lambda(Y)$, then $\Delta$ automatically scales
the measure $\eta$ and we denote the scaling constant by $\module \Delta$. Define
\begin{align*}
\Scentr := \{ H \subset \R \mid \; & \text{there exists $\; \Lambda
\actson (Y,\eta) \;$ free ergodic m.p.\ action,} \\ & \text{with
$\Lambda$ amenable and $\; \module(\Centr_\Lambda(Y)) = \exp(H)$} \;
\} \; .
\end{align*}

\begin{proposition}\label{prop.restrictions}
If $H \in \Sfactor$ or $H \in \Sequiv$, the subset $H \subset \R$ is a Borel
set that is Polishable: $H$ carries a unique Polish topology that is compatible
with its Borel structure.
\end{proposition}
\begin{proof}
The fundamental group
$\cF(M)$, resp.\ $\cF(\cR)$, of an arbitrary separable II$_1$ factor
$M$, resp.\ countable m.p.\ ergodic equivalence relation $\cR$ on a
standard probability space, appears in a short exact sequence
similar to \eqref{eq.exactsequence}. More precisely, denote by
$M^\infty := \B(\ell^2(\N)) \ovt M$ the amplified II$_\infty$
factor. The automorphism group $\Aut(M^\infty)$ is a Polish group
and consists of trace scaling automorphisms. The natural
homomorphism $\module : \Aut(M^\infty) \recht \Rp$ is continuous.
Therefore, the trace preserving automorphisms form a closed subgroup
of $\Aut(M^\infty)$ and $\cF(M)$ appears as the image of a Polish
group under an injective, continuous homomorphism. So, by a theorem
of Lusin and Souslin (see e.g.\ \cite[Theorem 15.1]{kechris}), the
group $\cF(M)$ is a Borel subset of $\Rp$. By construction, $\cF(M)$ is Polishable.
The same reasoning can be made for $\Sequiv$.
\end{proof}

\begin{remark} Proposition \ref{prop.restrictions} provides the only known a priori
restriction on the elements of $\Sfactor$ and $\Sequiv$, although we do not believe
that all Polishable Borel subgroups of $\R$ belong to $\Sfactor$ or $\Sequiv$.
On the other hand, more properties are known for the groups $H \in \bS$, but
the reason for this is rather indirect. For an arbitrary ergodic measure $\nu$,
a duality argument
(see \cite[Thm.~3.1]{aar-nad}) shows that the group $H_\nu$ also arises as the
\emph{eigenvalue group of a non-singular ergodic flow~:} there exists a non-singular,
ergodic action $(\sigma_t)_{t \in \R}$ on the standard infinite measure space $(Y,\eta)$
such that $H_\nu$ is exactly the group of $s \in \R$ for which there exists a non-zero
$F \in \rL^\infty(Y,\eta)$ satisfying $\sigma_t(F) = e^{its} F$ for all $t \in \R$.
Observe that the set of eigenvalue groups of non-singular ergodic flows, coincides
with the set
$$\bS_T = \{ H \subset \R \mid \; \text{there exists a factor $M$ with separable
predual such that}\; T(M) = H \} \; .$$ Here, $T(M)$ denotes Connes'
$T$-invariant of the factor $M$ (\cite{connesthesis}). So, $\bS
\subset \bS_T$, cf.\ \cite[Section~I]{giord-skand}.

Apart from being Polishable Borel subsets of $\R$, the following can be said about
the elements of $\bS_T$.
\begin{itemize}
\item Every $H \in \bS_T$ is a \emph{$\sigma$-compact} subset of $\R$.
This follows from the weak$^*$ compactness of the unit ball of $\rL^\infty(Y,\eta)$.
\item Every $H \in \bS_T$ is a \emph{saturated subgroup of $\R$.}
Following \cite[Section 2.1]{HMP}, a Borel set and subgroup $H \subset \R$
is called saturated if every bounded signed measure $\mu \in \operatorname{M}(\R)$
satisfies $|\mu(H)| \leq \sup \{|\widehat{\mu}(t)| \mid t \in \R \}$,
where $\widehat{\mu}$ denotes the Fourier transform of $\mu$.
\end{itemize}
Note in passing that the three properties of being Polishable,
$\sigma$-compact or saturated are rather independent.
Section 2.4 in \cite{HMP} provides $\sigma$-compact subgroups
of $\R$ that are not saturated. If $K \subset \R$ is a compact
\emph{Kronecker set\footnote{A compact subset $K \subset \R$
is called a Kronecker set if for every continuous function
$f : K \recht S^1$ and every $\eps > 0$, there exists
$y \in \R$ such that $|f(x) - \exp(ixy)|< \eps$ for all $x \in K$.}},
the subgroup generated by $K$ is a $\sigma$-compact subset of $\R$,
saturated (see \cite[Section 2.1]{HMP}), but not Polishable
by a Baire category argument (see e.g.\ \cite[Theorem 3]{LG}
or \cite[Remark 1.2]{aar-nad}).

All this makes it somewhat premature to set forth a plausible
conjecture on what should be an abstract characterization of the
groups in $\Sfactor$ and $\Sequiv$. Note however that all the groups
in these classes that we get in this paper are realized as
fundamental groups of II$_1$ factors and equivalence relations
arising from free ergodic m.p.\ actions of $\F_\infty$. Denoting by
$\Sfactor(\F_\infty)$, $\Sequiv(\F_\infty)$ the set of subgroups $H
\subset \R$ for which there exists a free ergodic m.p.\ action
$\F_\infty \actson X$ such that $\cF(\rL^\infty(X)\rtimes
\F_\infty)=\exp(H)$, respectively $\cF(\cR(\F_\infty\actson
X))=\exp(H)$, it seems very likely that one actually has
$\Sfactor=\Sequiv=\Sfactor(\F_\infty)=\Sequiv(\F_\infty)$. This
would provide a new universality property of $\F_\infty$, to be
compared with \cite{PS}.

In exactly the same way, we associate the sets $\Sfactor(\Gamma)$,
$\Sequiv(\Gamma)$ to any given countable group $\Gamma$. These
invariants capture interesting complexity aspects of the
group $\Gamma$, which we will investigate in a future paper. For
now, let us point out that by \cite{gab2}, if $\Gamma$ has at least
one $\ell^2$-Betti number not equal to $0$ or $\infty$ then
$\Sequiv(\Gamma)$ consists of the trivial group $\{1\}$ only. In
particular, $\Sequiv(\F_n) = \{\{1\}\}$, for any finite $n\geq 2$.
Also, one has $\{1\} \in \Sfactor(\F_n)$ by \cite{P-Betti}. We
expect that in fact $\Sfactor(\F_n)=\{\{1\}\}$ as well. On the other
hand, it is shown in \cite{P-corr}, \cite{NPS}, by using
``separability arguments'' inspired by Connes' pioneering rigidity
result in  \cite{connes}, that if $\Gamma$ is an infinite conjugacy
class group with the property (T), then $\Sfactor(\Gamma)$,
$\Sequiv(\Gamma)$ consist of countable groups only.

\end{remark}

\subsection*{Ergodic measures and Hausdorff dimension}

Examples of groups in $\bS$ can be constructed as follows (see e.g.\ \cite{aar-nad}).
Fix a sequence $a_n \in \N \setminus \{0,1\}$. Every $x \in \R$ can be uniquely written as
$$x = x_0 + \sum_{n=1}^\infty \frac{x_n}{a_1 \cdots a_n} \quad\text{where}\;\;
x_0 \in \Z, x_n \in \{0,\cdots,a_n-1\} \;\text{for all}\; n \geq 1 \; ,$$
and where we follow the convention that the sequence $(x_n)$ is not eventually
given by $x_n = a_n-1$.

Whenever $K_n \subset \{0,\cdots,a_n - 1\}$, define the generalized Cantor set
$$\cW_K := \{ x \in [0,1) \mid x_n \in K_n \;\text{for all}\; n \geq 1 \}$$
and equip $\cW_K$ with the probability measure $\lambda_K$ arising by identifying
(up to a countable set) $\cW_K$ and $\prod_{n=1}^\infty K_n$, the latter being
equipped with the product of the uniform probability measures.

Set $\al_n = a_1 \cdots a_n$.
Define the countable dense subgroup $Q \subset \R$ as the union of all $\al_n^{-1} \Z$.
It is straightforward to check that whenever $x \in Q$ and $\cU \subset \cW_K$ is
a Borel set such that $x + \cU \subset \cW_K$, then $\lambda_K(x + \cU) = \lambda_K(\cU)$.
Therefore, there is a unique $\sigma$-finite measure $\nu_K$ on the Borel sets of $\R$
that is $Q$-invariant, supported on $\cW_K + Q$ and such that the restriction of $\nu_K$
to $\cW_K$ equals $\lambda_K$. Moreover, the restriction of the $Q$-orbit equivalence
relation to $\cW_K \subset \R$ transfers to the equivalence relation on $\prod_{n=1}^\infty K_n$
given by
$$(x_n) \sim (y_n) \;\;\text{if and only if there exists $n_0$ such that}\;\; x_n =
y_n \;\;\text{for all}\;\; n \geq n_0 \; .$$
It follows that $\nu_K$ is an ergodic measure on $\R$.

Under certain conditions, the group $H_{\nu_K}$ can be determined explicitly. For example,
it follows from Theorem 4.1 in \cite{aar-nad} that if we take $K_n = \{0,\ldots,b_n -1\}$
such that $b_n < \frac{1}{2} a_n$ for $n$ large enough and $\sum_{n=1}^\infty b_n^{-1}
< \infty$, we have
\begin{equation}\label{eq.groupHnu}
H_{\nu_K} = \bigl\{x \in \R \;\big| \; \sum_{n=1}^\infty \frac{a_n}{b_n} \langle
\al_{n-1} x \rangle < \infty \bigr\}
\end{equation}
where $\langle x \rangle \in [0,1)$ denotes the distance of the real number $x$ to $\Z$.

If moreover $b_n \sim C_1 \rho_1^n$ and $a_n \sim C_2 \rho_2^n$ with $0 < C_1,C_2$
and $1 < \rho_1 < \rho_2$, one can prove as follows that the Hausdorff dimension
of $H_{\nu_K}$ is given by $\log \rho_1 / \log \rho_2$. First of all, by
\cite[Ch.~II, \S 7, Thm.~V]{kahane-salem}, the Hausdorff dimension of $\cW_K$
equals $\log \rho_1 / \log \rho_2$. Defining $$K'_n = \{0,\ldots,b_n - 1\} \cup
\{a_n - b_n,\ldots,a_n -1\} \; ,$$ it is clear that for every $x \in H_{\nu_K}$,
we eventually have $x_n \in K'_n$. So, $H_{\nu_K} \subset Q + \cW_{K'_n}$ and the
Hausdorff dimension of $H_{\nu_K}$ is at most $\log \rho_1 / \log \rho_2$. On the
other hand, for every $\rho_1^{-1} < \gamma < 1$, define $b\dpr_n$ as the integer
part of $\gamma^n b_n$. So, $b\dpr_n \sim C_1 (\gamma \rho_1)^n$. Set $K\dpr_n =
\{0,\ldots,b\dpr_n - 1 \}$. If $x \in \cW_{K\dpr}$, we have $$\langle \al_{n-1} x
\rangle \leq \frac{\gamma^n b_n}{a_n}$$
and so, $\cW_{K\dpr} \subset H_{\nu_K}$. It follows that the Hausdorff dimension
of $H_{\nu_K}$ is at least $\log(\gamma \rho_1) / \log(\rho_2)$ and this for all
$\rho_1^{-1} < \gamma < 1$. So, we have proven that the Hausdorff dimension of
$H_{\nu_K}$ is exactly $\log \rho_1 / \log \rho_2$.

\subsection*{Rigid actions}

We recall from 4.1 in \cite{P-Betti} the definition of relative
property (T) (or rigidity) for an inclusion of finite von Neumann
algebras and actions of groups.

\begin{definition} \label{def.relative-T}
Let $M$ be a factor of type II$_1$ with normalized trace
$\tau$ and let $A \subset M$ be a von Neumann subalgebra.
The inclusion $A \subset M$ is called \emph{rigid} if the following
property holds: for every $\eps > 0$, there exists a finite subset
$\cJ \subset M$ and a $\delta > 0$ such that whenever $_M H_M$ is
a Hilbert $M$-$M$-bimodule admitting a unit vector $\xi$ with the properties
\begin{itemize}
\item $\| a \cdot \xi - \xi \cdot a \| < \delta$ for all $a \in \cJ$,
\item $| \langle \xi, a \cdot \xi \rangle - \tau(a)| < \delta$ and
$|\langle \xi , \xi \cdot a \rangle - \tau(a) | < \delta$ for all
$a$ in the unit ball of $M$,
\end{itemize}
there exists a vector $\xi_0 \in H$ satisfying $\|\xi - \xi_0\| < \eps$
and $a \cdot \xi_0 = \xi_0 \cdot a$ for all $a \in A$.
\end{definition}

A free ergodic m.p.\ action $\Lambda \actson (X,\mu)$ is called
\emph{rigid} if the corresponding inclusion $\rL^\infty(X) \subset
\rL^\infty(X) \rtimes \Lambda$ is rigid. Recall from Proposition 5.1
in \cite{P-Betti} that if a group $\Gamma$ acts outerly on a
discrete abelian group $H$ and $\Gamma \actson (\hat{H},\text{\rm
Haar})$ is the action it induces on the (dual) compact group
$\hat{H}$ with its Haar measure, then $\Gamma \actson \hat{H}$ has
the relative property (T) iff the pair of groups $H \subset \Gamma
\ltimes H$ has relative property (T) of Kazhdan-Margulis.

Note that, by combining Theorems 4.4 in \cite{P-Betti} and A.1 in
\cite{NPS}, it follows that for every rigid action $\Lambda \actson
(X,\mu)$, the group $\Out(\cR(\Lambda^\infty \actson Y))$ is
countable. The same arguments show that if $\Lambda \actson (X,\mu)$
is a free m.p.\ action whose restriction to $\Gamma < \Lambda$ is
ergodic and rigid, the set $\Embed(\Gamma^\infty,\Lambda^\infty)$ is
countable modulo left composition with elements of
$[\Lambda^\infty]$.

Let $\Gamma_0 < \SL_2(\Z)$ be a non-amenable subgroup. By
\cite[Example 2 on page 62]{burger}, the pair $\Z^2 \subset \Z^2 \rtimes \Gamma_0$
has the relative property (T) in the group sense. Therefore, the natural
action $\Gamma_0 \actson \T^2$ is rigid.

\subsection*{Making actions freely independent}

A crucial role in our construction in the next section of an action
of $\F_\infty$ with ``special properties'' is played by choosing
inductively the action of the $n$-th generator of $\F_\infty$ from a
certain $G_\delta$-dense subset of $\Aut(X,\mu)$, the existence of
which uses a result in \cite{IPP}, \cite{To}, stated as Theorem
\ref{thm.makefree} below, for convenience (see A.1 in \cite{IPP} and
\cite{To} for details).

Let $(X,\mu)$ be the standard non-atomic probability space. We call
\emph{non-singular partial automorphism} of $(X,\mu)$ any
non-singular isomorphism $\phi$ between measurable subsets $D(\phi)
= X_0$ and $X_1$ of $X$. We identify partial automorphisms that are
equal almost everywhere. The partial automorphisms of $(X,\mu)$ form
a pseudo-group under composition. A partial automorphism $\paut$
with domain $D(\phi) \subset X$ is called (essentially) {\it free}
if the set $\{x \in D(\phi) \mid \paut(x) = x \}$ has measure zero.

\begin{theorem} \label{thm.makefree}
Let $\cV$ be a countable set of free non-singular partial
automorphisms of $(X,\mu)$. Define $T \subset \Aut(X,\mu)$ as the
set of measure preserving automorphisms $\theta$ such that every
non-empty word with letters alternatingly from $\cV$ and $\theta \cV
\theta^{-1}$ defines a free partial automorphism of $(X,\mu)$. Then,
$T$ is a $G_\delta$-dense subset of the Polish group of m.p.\ automorphisms of $(X,\mu)$.
\end{theorem}

We will use Theorem \ref{thm.makefree} under the following form. Set
$(Y,\eta) = (X \times \Z, \mu \times c)$. Suppose that $\cW$ is a
countable set of non-singular partial automorphisms of $Y$ with the
property that $\si_n \circ \paut$ is free for all $n \in \Z$ and
$\paut \in \cW$. Here $(\si_n)_{n \in \Z}$ denotes the shift action
$\Z \actson \Z$. Define $T \subset \Aut(X,\mu)$ as the set of
measure preserving automorphisms $\theta$ of $(X,\mu)$ such that
every non-empty word with letters alternatingly from $\cW$ and
$\theta \cW \theta^{-1}$ defines a free partial automorphism of $Y$.
Then, $T$ is a $G_\delta$-dense subset of $\Aut(X,\mu)$. This
statement follows by defining the partial isomorphisms $\paut_n : X
\recht X \times \{n\} \subset Y : \paut_n(x) = (x,n)$ and applying
Theorem \ref{thm.makefree} to the countable set $\cV = \{
\paut_n^{-1} \paut \paut_m \mid n,m \in \Z , \paut \in \cW \}$ of
free partial automorphisms of $(X,\mu)$.

\section{Construction of the basic $\F_\infty$-action}

The actions $\F_\infty\actson X_i$ in the statement of Theorem
1.1 will (roughly) be obtained as diagonal product actions of $\F_\infty$, the first
component of which will be referred to as the
{\it basic action}. Its construction, which is the subject of this
section, is inspired by arguments in \cite{GP}, \cite{IPP}.

To start with, let $H = \F_\infty$ with free generators
$a,b,g_1,g_2,\ldots$. Define as follows subgroups of $H$. We set
$H_0 := \langle a,b \rangle$ and, for every subset $E \subset \No$,
$H_E := \langle a, b, g_n \mid n \in E \rangle$. We set $H_n =
\langle a,b, g_1,\ldots, g_n \rangle$.

Define $(X,\mu)$ to be $\T^2$ with the Lebesgue measure and
$H_0 \actson X$ by viewing $H_0 \subset \SL_2(\Z)$
and restricting to $H_0$ the canonical action of $\SL(2,\Z)$
on $\T^2$. Note that the action $H_0 \actson X$ is rigid and
weakly mixing (see e.g.\ \cite[Corollary 3.3.2]{P-Betti}). Denote as before
$(Y,\eta) = (X \times \Z, \mu \times c)$,
where $c$ is the counting measure on $\Z$.

Let $K$ be any countable group and $K \actson (X,\mu)$ an arbitrary
free m.p.\ action. Use Theorem \ref{thm.makefree} to conjugate $K
\actson (X,\mu)$ such that the action $H_0 * K \actson (X,\mu)$ is
still free.

We construct, by induction on $n$, a free m.p.\ action of $H_n * K$
on $(X,\mu)$, extending the action of $H_0 * K$ chosen above.
Suppose that we are given a free m.p.\ action of $H_n * K$ on $X$
extending the $H_0 * K$-action chosen above. Choose for every $F
\subset \{1,\ldots,n\}$ a set $\cA(F) \subset \Emb$ of
representatives for the outer elements in $\Emb$ modulo
$[(H_F*K)^\infty]$. Note that $\cA(F)$ is countable, because $H_F
\actson X$ is a rigid action. Consider the following countable set
of non-singular automorphisms of $Y$~:
\begin{align*}
\cV := & \{\si_g \Delta^{\pm 1} \si_h \mid g,h \in (H_n*K)^\infty \; ,
\; \Delta \in \cA(F) \;\;\text{for some}\;\; F \subset \{1,\ldots,n\} \} \\
& \cup \{ \si_k \mid  k \in (H_n * K - \{e\}) \times \Z \} \; .
\end{align*}

\begin{lemma}
The automorphism $\si_n \circ \Delta$ is free whenever $n \in \Z$
and $\Delta \in \cV$.
\end{lemma}
\begin{proof}
The lemma follows once we show that $\si_{g^{-1}} \circ \Delta$ is
free whenever $g \in (H_n*K)^\infty$ and $\Delta \in \Emb -
[(H_F*K)^\infty]$. If this is not the case, we find $\Delta(x) = g
\cdot x$ for all $x$ in some non-negligible subset $\cU \subset Y$.
Because $H_F^\infty \actson Y$ is ergodic and $\Delta$ belongs to
$\Emb$, it follows that $\Delta(x) \in (H_n*K)^\infty \cdot x$ for
almost all $x \in Y$. We have to prove that actually $\Delta(x) \in
(H_F*K)^\infty \cdot x$ almost everywhere, in order to reach a
contradiction with the outerness of $\Delta$. Define $\Delta_1 : X
\recht X$ such that $\Delta(x,0) \in \{\Delta_1(x)\} \times \Z$ for
almost all $x \in X$. Since $\Delta(x) \in (H_n*K)^\infty \cdot x$
for almost all $x \in Y$, we find $\vphi : X \recht H_n * K$ such
that $\Delta_1(x) = \vphi(x) \cdot x$ for almost all $x \in X$.
Since $\Delta \in \Emb$, it follows that $\vphi(g \cdot x) g
\vphi(x)^{-1} \in H_F * K$ for almost all $x \in X$, $g \in H_F$. By
weak mixing of $H_F \actson (X,\mu)$ and because an element $h \in
H_n * K$ satisfying $h H_F h^{-1} \subset H_F * K$ must belong to
$H_F * K$, it follows that $\vphi(x) \in H_F * K$ for almost all $x
\in X$. So, we are done.
\end{proof}

Combining Theorem \ref{thm.makefree}, the remarks following that
theorem and the previous lemma, take a free ergodic m.p.\ action of
$\Z \cong g_{n+1}^\Z$ on $(X,\mu)$ such that every non-empty word
with letters alternatingly from $\cV$ and $\{\si_{g_{n+1}^k} \mid k
\in \Z - \{0\}\}$ yields a free transformation of $Y$.

We obtained in particular a free m.p.\ action of $H_{n+1} * K$
extending the $H_n * K$-action that we started with. Continuing by
induction, we obtain a  free m.p.\ action $\si : H * K \actson
(X,\mu)$, whose restriction to $H_0$ is ergodic and rigid.

\section{Actions with no symmetries}

The following theorem is the crucial technical result on which all results
in the paper rely.

As above $H \cong \F_\infty$, $K$ is arbitrary and we defined the action
$\si : H * K \actson (X,\mu)$. We defined the subgroup
$H_E \subset H$ generated by $a,b,g_n, n \in E$ whenever $E \subset \No$.
All $H_E$ act rigidly on $(X,\mu)$.

\begin{theorem} \label{thm.main-technical}
There exists an uncountable family $\cE$ of infinite subsets of $\N$
such that $$\Emb[E] = [(H_E * K)^\infty]$$ for all $E \in \cE$.
In particular, $\Out(\cR((H_E*K)^\infty \actson Y))$ is trivial for all $E \in \cE$.
\end{theorem}

\begin{proof}
Choose an uncountable family $\cE_1$ of infinite subsets of $\N$ satisfying
$E \cap F$ finite whenever $E,F \in \cE_1$ and $E \neq F$. Suppose that the theorem
does not hold. Leaving out the countably many $E \in \cE_1$ such that
$\Emb[E] = [(H_E * K)^\infty]$, we find an
uncountable subset $\cE \subset \cE_1$ and for every $E \in \cE$ an outer
$\Delta_E \in \Emb[E]$.

Recall that we fixed, for every finite subset $F \subset \N$, a set of
representatives $$\cA(F) \subset \Emb$$ for the outer elements in
$\Emb$ modulo $[(H_F * K)^\infty]$.

{\bf Step 1.} {\it There exist $E,F \in \cE$ such that $E \neq F$ and
$\Delta_E(x) = \Delta_F(x)$ for all $x$ in a non-negligible subset
$\cU \subset Y$.} Step~1 will follow by using relative property~(T)
and a separability argument. Denote $A = \rL^\infty(X,\mu)$ and define
the group measure space II$_1$ factors $M_0 = A \rtimes H_0$ and
$M = A \rtimes (H*K)$. We write $M^\infty := \B(\ell^2(\Z)) \ovt M$
and view $\rL^\infty(Y) = A^\infty$ as the obvious diagonal subalgebra
of $M^\infty$. Denote by $\Tr$ the natural semi-finite trace on $M^\infty$.
For every $E \in \cE$, define $p_E \in A^\infty$ as the projection on
$\Delta_E(X \times \{0\})$. Note that $\Tr(p_E) < \infty$. Define the embeddings $\beta_E : M_0 \recht p_E M^\infty p_E$ by naturally
extending the homomorphism $\beta_E : A \recht A^\infty p_E : \beta_E(a) =
a \circ \Delta_E^{-1}$ for all $a \in A$. Define then for all $E,F \in \cE$
the Hilbert $M_0$-$M_0$-bimodule $\cH_{E,F} := p_E \rL^2(M^\infty) p_F$ with
bimodule action given by
$$a \cdot \xi \cdot b := \beta_E(a) \xi \beta_F(b) \; .$$
Set $\xi_{E,F} = \Tr(p_E p_F)^{-1/2} p_E p_F$ whenever $p_E p_F$ is non-zero.
As mentioned above, the inclusion $A \subset M_0$ is rigid. Set $\eps = 1/2$
and take the finite subset $\cJ \subset M_0$ and $\delta > 0$ such that the
conclusion of Definition \ref{def.relative-T} holds. The separability of the
Hilbert space $\rL^2(M^\infty)$ and the uncountability of $\cE$ provide us
with $E,F \in \cE$, $E \neq F$ such that the vector $\xi_{E,F}$ satisfies
the assumptions in Definition \ref{def.relative-T}. So, there exists a vector
$\xi_0 \in \cH_{E,F}$ satisfying $a \cdot \xi_0 = \xi_0 \cdot a$ for all $a \in A$
and $\|\xi_{E,F} - \xi_0 \| < \eps = 1/2$. Taking the orthogonal projection of
$\xi_0$ onto $\rL^2(Y)$, we may assume that $\xi_0$ has its matrix coefficients in $A$.
Since $\xi_0$ is non-zero, the claim of Step~1 follows.

{\bf Step 2.} {\it There exists $\paut \in [(H_E*K)^\infty]$ such that
$\paut \circ \Delta_E \in \cA(E \cap F)$.}
We have $\Delta_E(x) = \Delta_F(x)$ for all $x \in \cU$. Since $H_0^\infty
\actson Y$ is ergodic and $H_0 \subset H_E \cap H_F$,
it follows that $$\Delta_E(x) \in \bigl((H_E*K)^\infty \, (H_F * K)^\infty \bigr) \,
\cdot
\Delta_F(x)$$ for almost all $x \in Y$. Composing $\Delta_E$ and $\Delta_F$
on the left by elements in resp.\ $[(H_E*K)^\infty]$ and $[(H_F*K)^\infty]$,
we may assume that $\Delta_E(x) = \Delta_F(x)$ for almost all $x \in Y$. But
then, $\Delta_E \in \Emb[E \cap F]$. Since $\Delta_E$ is outer, we can compose
$\Delta_E$
by an element of $[(H_{E \cap F}*K)^\infty]$ and assume that
$\Delta_E \in \cA(E \cap F)$. This concludes Step~3.

{\bf Step 3.} {\it We obtain a contradiction with our choice of action
$\si : H * K \actson X$.} Since $E$ is infinite and $E \cap F$ finite,
take $n \in E$ such that $k < n$ for all $k \in E \cap F$. Since
$\Delta_E$ belongs to $\Emb[E]$,
we can take $h \in (H_E*K)^\infty$ and a non-negligible subset
$Y_0 \subset Y$ such that $\Delta_E(g_n \cdot x) = h \cdot \Delta_E(x)$
for all $x \in Y_0$. Suppose first that $h \in (H_n * K)^\infty$.
Write $h=(h_0,m)$ where $m \in \Z$ and $h_0$ is a word with letters
alternatingly in $H_{n-1} * K - \{e\}$ and $\{ g_n^k \mid k \in \Z - \{0\} \}$. The fact that
$$\si_{h^{-1}} \circ \Delta_E \circ \si_{g_n} \circ \Delta_E^{-1}$$
is not free is a contradiction with $\Delta_E$ belonging to $\cA(E
\cap F)$ and our choice of $\si_{g_n}$. Suppose next that $h \not\in
(H_n * K)^\infty$ and take the smallest $p$ such that $h \in (H_p *
K)^\infty$. So, $p > n$ and the generator $g_p$ appears in the
reduced expression of $h$. Moreover,
$$\si_{g_n} \circ \Delta_E^{-1} \circ \si_{h^{-1}} \circ \Delta_E$$
is not free. This is also a contradiction with our choice of
$\si_{g_p}$, because $\Delta_E$ belongs to $\cA(E \cap F)$.
\end{proof}

Denote by $\si_E$ the restriction of $\si$ to $H_E * K$.

\begin{lemma} \label{lem.non-oe}
Let $E \subset \N$. There are at most countably many $F \subset \N$
such that $\si_E$ and $\si_F$ are stably orbit equivalent.
\end{lemma}
\begin{proof}
Let $E \subset \N$. Suppose that $\cE$ is an uncountable family of subsets of $\N$ such
that $\si_E$ and $\si_F$ are stably orbit equivalent for all $F \in \cE$.
Take for every $F \in \cE$, a non-singular automorphism $\Delta_F : Y \recht Y$
satisfying $\Delta_F( (H_E * K) \cdot x) = (H_F * K) \cdot \Delta_F(x)$ for
almost all $x \in Y$. Define the II$_1$ factors $M_0$ and $M$ as in Step~1
of the proof of \ref{thm.main-technical}. The formula $\beta_F(a) = a \circ \Delta_F^{-1}$
defines, as in the proof of \ref{thm.main-technical}, a homomorphism $\beta_F : A \recht
A^\infty p_F$ that naturally extends to a homomorphism $\beta_F : M_0 \recht
p_F M^\infty p_F$. Here $A^\infty = \rL^\infty(Y)$ and $p_F$ is the projection
onto $\Delta_F(X \times \{0\})$. Separability and rigidity provide us with
$F, F' \in \cE$, $F \neq F'$ such that $\Delta_F(x) = \Delta_{F'}(x)$ for
all $x$ belonging to a non-negligible subset $\cU \subset Y$. As in the
proof of \ref{thm.main-technical}, we can compose $\Delta_F,\Delta_{F'}$ on the
left with elements in $[(H_{F}*K)^\infty]$, $[(H_{F'}*K)^\infty]$ and
assume that $\Delta_F(x) = \Delta_{F'}(x)$ for almost all $x \in Y$. From
this it follows that $H_F*K$ and $H_{F'}*K$ are equal subgroups of $H * K$, a contradiction.
\end{proof}

Combining Theorem \ref{thm.main-technical} and Lemma \ref{lem.non-oe},
we immediately get the following result.

\begin{theorem} \label{thm.final}
Let $\si : K \actson (X,\mu)$ be any free, probability measure
preserving action of any countable group. There exists an
uncountable family $(\si_i)_{i \in I}$ of free ergodic m.p.\ actions
of $\F_\infty * K$ on probability spaces $(X_i,\mu_i)$ with the
following properties:
\begin{itemize}
\item The orbit equivalence relation of $\F_\infty * K \actson X_i$
has trivial fundamental group
and trivial outer automorphism group.
\item The restriction of $\si_i$ to $K$ is conjugate with $\si$.
\item The restriction of $\si_i$ to $\F_\infty$ is ergodic and rigid.
\item The actions $(\si_i)_{i \in I}$ are not stably orbit equivalent.
\end{itemize}
\end{theorem}

\section{Actions with prescribed symmetries}

We use the following terminology and notations. When $(Y,\eta)$ is a standard
infinite measure space and $\Lambda \actson (Y,\eta)$ is a measure preserving,
ergodic action, we define for every $\al \in \Aut(\Lambda)$,
the group $\Aut^\al_\Lambda(Y)$ of non-singular automorphisms $\Delta$ of $Y$
satisfying $\Delta(g \cdot y) = \al(g) \cdot \Delta(y)$ for all $g \in \Lambda$
and almost all $y \in Y$. We set $\Centr_\Lambda(Y) = \Aut^{\id}_\Lambda(Y)$.

We say that $\Lambda \actson (Y,\eta)$ is induced from $\Lambda_0 \actson Y_0$
if $Y_0$ is a measurable subset of $Y$ that is globally $\Lambda_0$-invariant
and satisfies $\eta(g \cdot Y_0 \cap Y_0) = 0$ whenever $g \not\in \Lambda_0$.

\begin{lemma} \label{lemma.main}
Suppose that we have the following data.

\begin{itemize}
\item A standard infinite measure space $(Y,\eta)$ with a  non-singular action
$\Lambdatil \actson (Y,\eta)$ scaling the measure $\eta$.
\item A normal subgroup $\Lambda \lhd \Lambdatil$ such that the restricted action
$\Lambda \actson (Y,\eta)$ is measure preserving, free and ergodic.
\end{itemize}
We make the following additional assumption.
\begin{itemize}
\item The action $\Lambda \actson (Y,\eta)$ defines the same orbit equivalence
relation on $(Y,\eta)$ as the free, ergodic, measure preserving
action $\Sigma \actson (Y,\eta)$ that moreover is profinite in the
following sense: we have $(Y,\eta) = \invlimit (Y_n,\eta_n)$ where
the $(Y_n,\eta_n)$ are atomic infinite measure spaces and the
subalgebras $\rL^\infty(Y_n,\eta_n) \subset \rL^\infty(Y,\eta)$ are
globally $\Sigma$-invariant.
\end{itemize}
Note that the additional assumption is automatic when $\Lambda$ is
amenable, but of course also when $\Lambda \actson (Y,\eta)$ is itself profinite.

Suppose next that $\Gammatil \actson (X,\mu)$ is a free m.p.\ action
of the group $\Gammatil$ on a probability measure space $(X,\mu)$
and $H$ is a subgroup of $\Gammatil$ with the following properties:
\begin{itemize}
\item The restricted action $H \actson (X,\mu)$ is ergodic and rigid.
\item We have $\Embed(H^\infty,{\Gammatil}^\infty) = [{\Gammatil^\infty}]$.
\end{itemize}

Whenever $\pi : \Gammatil \recht \Lambdatil$ is a surjective
homomorphism with $H \subset \Ker \pi$, we set $\Gamma =
\pi^{-1}(\Lambda)$ and define $\cR_\pi$ as the equivalence relation
given by the orbits of the free, ergodic, (infinite) measure
preserving action $\Gamma \actson X \times Y : g \cdot (x,y) = (g
\cdot x , \pi(g) \cdot y)$.

The following results hold.
\begin{enumerate}
\item \label{part-one} We have $$\Out(\cR_\pi) \cong \frac{\{(g,\Delta) \mid g
\in \Lambdatil , \Delta \in \Aut^{\Ad g}_\Lambda(Y) \}}{\Lambda} \; .$$
In particular, we have a short exact sequence
    $$e \recht \Centr_\Lambda(Y) \recht \Out(\cR_\pi) \recht \frac{\Lambdatil}{\Lambda}
    \recht e \; .$$
\item \label{part-two} The equivalence relations $\cR_{\pi_1}$ and $\cR_{\pi_2}$
are stably orbit equivalent if and only if $\Lambda \actson Y$ is induced from
$\Lambda_i \actson Y_i \; (i=1,2)$ in such a way that $\pi_1^{-1}(\Lambda_1) =
\pi_2^{-1}(\Lambda_2)$ and there exists a non-singular isomorphism $\Delta_0 :
Y_1 \recht Y_2$ satisfying $\Delta_0(\pi_1(g) \cdot y) = \pi_2(g) \cdot \Delta_0(y)$
for all $g \in \pi_i^{-1}(\Lambda_i)$ and almost all $y \in Y_1$.

    In particular, if $\Lambda \actson Y$ can only be induced from faithful actions,
    stable orbit equivalence of $\cR_{\pi_1}$ and $\cR_{\pi_2}$ implies $\Ker \pi_1 =
    \Ker \pi_2$.
\end{enumerate}
\end{lemma}

\begin{proof}
Choose a surjective homomorphism $\pi : \Gammatil \recht \Lambdatil$ and define
$\cR_\pi$ as the orbit equivalence relation of $\Gamma \actson X \times Y$.
When $g \in \Gammatil$ and when $\Delta_0$ is a non-singular automorphism of
$(Y,\eta)$ scaling the measure $\eta$ and satisfying $\Delta_0(h \cdot y) =
(\pi(g)h\pi(g)^{-1}) \cdot \Delta_0(y)$ for all $h \in \Lambda$ and almost
all $y \in Y$, we define the automorphism $\Delta$ of $X \times Y$ as $\Delta(x,y) =
(g \cdot x, \Delta_0(y))$ and observe that $\Delta$ normalizes the action $\Gamma
\actson X \times Y$. In particular, $\Delta$ is an automorphism of $\cR_\pi$.
Clearly, $\Delta$ is inner if and only if $g \in \Gamma$ and $\Delta_0(y) =
\pi(g) \cdot y$ for almost all $y$.

To prove item~\ref{part-one}, we have to show that every automorphism of $\cR_\pi$
is of the form above modulo $[\Gamma]$. The main step of the proof consists
in using the rigidity of $\Gamma \actson X$ and the profiniteness of $\Sigma
\actson Y$ in order to prove that any automorphism $\Delta(x,y) =
(\Delta_1(x,y),\Delta_2(x,y))$ of $\cR_\pi$ is such that $\Delta_1$
is essentially independent of the $y$-variable (cf.\ Section~6 in \cite{P-Betti}).

Set $A = \rL^\infty(X,\mu)$ and $B = \rL^\infty(Y,\eta)$. Denote $M = (A \ovt B)
\rtimes \Gamma$ and note that $M$ is a II$_\infty$ factor. Write $\cM = (A \rtimes \Gamma)
\ovt (B \rtimes \Lambda)$. The map $$(a \ot b)u_g \mapsto (au_g) \ot (b u_{\pi(g)})$$
extends to a trace preserving embedding $M \subset \cM$. Let $\Sigma \actson (Y,\eta)$
be a profinite action with $\Sigma \cdot y = \Lambda \cdot y$ for almost all $y \in Y$.
Identify $B \rtimes \Lambda$ with $B \rtimes \Sigma$. Take $(Y,\eta) = \invlimit (Y_n,\eta_n)$
in such a way that the subalgebras $B_n := \rL^\infty(Y_n,\eta_n) \subset B$ are globally
$\Sigma$-invariant. Denote by $E_n : B \rtimes \Sigma \recht B_n \rtimes \Sigma$ the unique
trace preserving conditional expectation.

Let $\Delta$ be an automorphism of $\cR_\pi$. Write $\Delta(x,y) =
(\Delta_1(x,y),\Delta_2(x,y))$. Let $p$ be a minimal projection in
$B_0 = \rL^\infty(Y_0,\eta_0)$ and rescale $\eta$ such that
$p$ has trace $1$. So, $p$ projects onto a measurable subset $Z \subset Y$ with
$\eta(Z) = 1$. Denote by $p_1$ the projection onto $\Delta(X \times Z)$. Note that $p$
and $p_1$ are finite projections
in the II$_\infty$ factor $M$. The automorphism $\Delta$ yields a $^*$-isomorphism of
II$_1$ factors
$$\Psi : p_1 M p_1 \recht p M p \quad\text{satisfying}\;\; \Psi(p_1(A \ovt B)) =
p (A \ovt B) \;\;\text{and}\;\; \Psi(f)(x,y) = f(\Delta(x,y))$$
for almost all $(x,y) \in X \times Z$ whenever $f \in p_1(A \ovt B)$.

Set $D = \psi(p_1(A \ot 1))$. Then, $D \subset pMp$ is a rigid inclusion and hence,
also the inclusion $D \subset p\cM p$ is rigid. Since $\|(\id \ot E_n)(a) - a \|_2 \recht 0$
for all $a \in p \cM p$, it follows that $\|(\id \ot E_n)(d) - d \|_2 \recht 0$ uniformly on
the unit ball of $D$. Lemma \ref{lem.intertwine} below allows to take $n$ and a projection $q$
onto a subset
$\cU \subset X \times Z$ of measure strictly greater than $\frac{3}{4}$ such
that $D q \subset (A \ovt B_n)q$. Denote by $\rho_n : Y \recht Y_n$ the quotient
maps given by $(Y,\eta) = \invlimit (Y_n,\eta_n)$. We get measurable subsets
$\cU_1 \subset X \times Y_n$ and $X_1 \subset X$
as well as a quotient map of measure spaces $\al : \cU_1 \recht X_1$ such that
$\Delta_1(x,y) = \al(x,\rho_n(y))$ for almost all $(x,y) \in \cU$.

Since also $(A \ovt B_n)p \subset pMp$ is rigid, we use $\Psi^{-1}$ and similarly
find for every $\eps > 0$, an $m \in \N$ and a projection $q'$ in $(A \ovt B)p$
with $\tau(q') \geq 1 - \eps$, such that
$$(A \ovt B_n) q' \subset \Psi(p_1 (A \ovt B_m)) q' \; .$$
This means that, after making $\cU$ slightly smaller but still with measure
greater than $\frac{3}{4}$, the quotient map $\al : \cU_1 \recht X_1$ may be
assumed to be countable-to-one.

Partitioning $Y$ in subsets on which $\rho_n$ is constant, we find a
non-negligible subset $Z_0 \subset Z$ and a measurable subset $\cV \subset
X \times Z_0$ with $(\mu \times \eta)(\cV) \geq \frac{3}{4} \eta(Z_0)$ such
that $\Delta_1(x,y) = \al(x)$ for almost all $(x,y) \in \cV$, where
$\al : X_2 \subset X \recht X_1 \recht X$ is a countable-to-one quotient map.
Whenever $x \in X$, define $\cV_x = \{y \in Z_0 \mid (x,y) \in \cV \}$. Let
$\cW \subset X$ be the necessarily non-negligible subset of $x \in X$ with
$\eta(\cV_x) > \frac{2}{3} \eta(Z_0)$. Whenever $g \in H$, $x \in \cW$ and
$g \cdot x \in \cW$, the set $\cV_{g \cdot x} \cap \cV_x$ has measure at
least $\frac{1}{3} \eta(Z_0)$. So, we can take $y \in \cV_{g \cdot x}
\cap \cV_x$ and hence
$$\al(g \cdot x) = \Delta_1(g \cdot x, y) = \Delta_1(g \cdot (x,y)) \in
\Gamma \cdot \Delta_1(x,y) = \Gamma \cdot \al(x) \; .$$

Since $\al$ is countable-to-one, we can further restrict $\al|_\cW$ to
a non-negligible partial automorphism of $(X,\mu)$ with the property
that $\al(g \cdot x) \in \Gamma \cdot \al(x)$ whenever $g \in H$ and
both $x$ and $g \cdot x$ belong to the domain of $\al$. By our assumption
on the action $\Gammatil \actson (X,\mu)$, it follows that $\al(x) \in
\Gammatil \cdot x$ for all $x$ in the domain of $\al$. So, we can compose
the automorphism $\Delta$ by an automorphism of the form $(x,y) \mapsto
(g \cdot x,\pi(g) \cdot y)$ where $g \in \Gammatil$ and we may assume
now that $\Delta(x,y) \in \{x\} \times Y$ for all $(x,y)$ in some non-negligible
subset $\cU$ of $X \times Y$. We prove that, module $[\Gamma]$, such a $\Delta$
has the form $\Delta(x,y) = (x,\Delta_0(y))$, where $\Delta_0 \in \Centr_\Lambda(Y)$.

Since the action $\Gamma \actson X \times Y$ is ergodic, choose a measurable map
$\vphi : X \times Y \recht \Gamma$ such that $\vphi(x,y) \cdot (x,y) \in \cU$ for
almost all $(x,y) \in X \times Y$. Define the measurable map
$$\psi : X \times Y \recht X \times Y : \psi(x,y) = \vphi(x,y)^{-1} \cdot
\Delta(\vphi(x,y) \cdot (x,y)) \; .$$
Note that $\psi(x,y) \in \{x\} \times Y$ for almost all $(x,y)$ and that $\psi$
preserves $\cR_\pi$ in the sense that $(x,y)\cR_\pi (x',y')$ if and only if
$\psi(x,y) \cR_\pi \psi(x',y')$. It follows that $\psi(g \cdot (x,y)) =
g \cdot \psi(x,y)$ for all $g \in \Gamma$ and almost all $(x,y)$. Also
note that $X \times Y$ can be partitioned into subsets $X \times Y =
\bigsqcup_m Z_m$ such that the restriction of $\psi$ to $Z_m$ is a non-singular
partial automorphism of $X \times Y$.

We claim that $\psi$ is a non-singular automorphism of $X \times Y$. First assume
that $\psi$ is not essentially injective. We get a partial automorphism $\theta$
of $X \times Y$ with domain and range being disjoint subsets of $X \times Y$ and
$\psi(\theta(x,y)) = \psi(x,y)$ for all $(x,y) \in D(\theta)$. But then,
$\theta(x,y) \cR_\pi (x,y)$ for almost all $(x,y) \in D(\theta)$. Making
the domain of $\theta$ smaller, we find $g \in \Gamma - \{e\}$ and a non-negligible
subset $\cV \subset X \times Y$ such that $\psi(g \cdot (x,y)) = \psi(x,y)$ for all
$(x,y) \in \cV$. Since $\psi(g \cdot (x,y)) = g \cdot \psi(x,y)$, we reached a
contradiction with the essential freeness of $\Gamma \actson X$. Since the range of
$\psi$ is $\Gamma$-invariant, $\psi$ is also essentially surjective.

It follows that $\psi(x,y) = (x, \psi_x(y))$, where $\psi_x$ is, for almost every
$x \in X$, a non-singular automorphism of $(Y,\eta)$. Since $\psi(g \cdot (x,y)) =
g \cdot \psi(x,y)$, it follows that $\psi_{g \cdot x}(\pi(g) \cdot y) = \pi(g) \cdot
\psi_x(y)$. In particular, $\psi_{g \cdot x} = \psi_x$ for all $g \in \Ker \pi$ and
almost all $x \in X$. Since $\Ker \pi$ acts ergodically on $(X,\mu)$, it follows that
almost every $\psi_x$ equals almost everywhere $\Delta_0$, where $\Delta_0 \in
\Centr_\Lambda(Y)$. This concludes the proof of item~\ref{part-one}.

We now prove item~\ref{part-two}. Let $\pi_1,\pi_2 : \Gammatil \recht \Lambdatil$
be surjective homomorphisms such that $H \subset \Ker \pi_i$ for $i=1,2$.

Suppose first that $\Lambda \actson Y$ is induced from $\Lambda_i \actson Y_i$ in
such a way that the conditions in item~\ref{part-two} hold. Write $\Gamma_0 =
\pi_i^{-1}(\Lambda_i)$. For $i=1,2$, we have actions $$\si_i : \Gamma_0 \actson X
\times Y_i : g \cdot (x,y) = (g \cdot x, \pi_i(g) \cdot y) \; .$$
The map $(x,y) \mapsto (x,\Delta_0(y))$ conjugates $\si_1$ and $\si_2$. Moreover,
the reduction of $\cR_{\pi_i}$ to $X \times Y_i$ is precisely given by the orbit
equivalence relation of $\si_i$. So, these reductions are orbit equivalent
and hence, $\cR_{\pi_1}$ is stably orbit equivalent with $\cR_{\pi_2}$.

Suppose next that $\Delta$ is a non-singular automorphism of $X \times Y$
defining a stable orbit equivalence of $\cR_{\pi_1}$ and $\cR_{\pi_2}$.
Write $\Gamma_i = \pi_i^{-1}(\Lambda)$. The same rigidity vs.\ profiniteness argument as above tells us that we can compose $\Delta$ with an automorphism of
$\cR_{\pi_2}$ and assume that $\Delta(x,y) \in \{x\} \times Y$ for all $(x,y)$
in some non-negligible subset $\cU$ of $X \times Y$.
Note that writing $\cV = \Delta(\cU)$, we also have $\Delta^{-1}(x,y) \in \{x\} \times Y$
for all $(x,y)$ in $\cV$.

Let $\cU_0$ be the set of $y \in Y$ such that $\cU^y := \{x \in X \mid (x,y) \in \cU\}$
is non-negligible. Define $\Lambda_1$ as the subgroup of $\Lambda$ generated by $g \in \Lambda$
satisfying
$\eta(g \cdot \cU_0 \cap \cU_0) > 0$. Set $Y_1 = \bigcup_{g \in \Lambda_1} g \cdot \cU_0$.
By construction, $\Lambda \actson Y$ is induced from $\Lambda_1 \actson Y_1$. Therefore,
the latter is ergodic and also the action $\pi_1^{-1}(\Lambda_1) \actson X
\times Y_1$ follows ergodic.

We also define $\cV_0$ as the set of $y \in Y$ such that $\cV^y$ is non-negligible.
This leads to a similar construction of $\Lambda_2 \actson Y_2$ from which $\Lambda \actson Y$
is induced.

We claim that $\pi_1^{-1}(\Lambda_1) = \pi_2^{-1}(\Lambda_2)$. By symmetry, it is
sufficient to take $g \in \pi_1^{-1}(\Lambda_1)$ and to prove that $\pi_2(g) \in
\Lambda_2$. It is even sufficient to start with $g \in \Gamma_1$ such that
$\pi_1(g)^{-1} \cdot \cU_0 \cap \cU_0$ is non-negligible. This means that
$(gH)^{-1} \overset{\si_1}{\cdot} \cU \cap \cU$ is non-negligible. Since $\Delta$
maps $\cR_{\pi_1}$ onto $\cR_{\pi_2}$ and satisfies $\Delta(x,y) \in \{x\}
\times Y$ for all $(x,y) \in \cU$, it follows that $g \in \Gamma_2$ and that
$(gH)^{-1} \overset{\si_2}{\cdot} \cV \cap \cV$ is non-negligible. Hence,
$\pi_2(g)^{-1} \cdot \cV_0 \cap \cV_0$ is non-negligible,
proving that $\pi_2(g) \in \Lambda_2$.

We set $\Gamma_0 = \pi_1^{-1}(\Lambda_1) = \pi_2^{-1}(\Lambda_2)$. By ergodicity of
$\Gamma_0 \actson X \times Y_1$, we can choose a
measurable map $\vphi : X \times Y_1 \recht \Gamma_0$ such that $\vphi(x,y)
\cdot (x,y) \in \cU$ for almost all $(x,y) \in X \times Y_1$. Write
$$\psi : X \times Y_1 \recht X \times Y_2 : \psi(x,y) = \vphi(x,y)^{-1} \cdot
\Delta(\vphi(x,y) \cdot (x,y)) \; .$$
In exactly the same way as in the proof of item~\ref{part-one}, it follows that
$\psi(x,y) = (x,\Delta_0(y))$,
where $\Delta_0$ is a non-singular isomorphism of $Y_1$ onto $Y_2$ satisfying $\Delta_0(\pi_1(g)
\cdot y) = \pi_2(g) \cdot \Delta_0(y)$ for all $g \in \Gamma_0$ and almost all $y \in Y$.
\end{proof}

As we will see below, Lemma \ref{lemma.main} allows us to construct
free, ergodic, measure preserving actions of groups on
\emph{infinite} measure spaces, whose orbit equivalence relation has
prescribed fundamental group and, to a certain extent, prescribed
outer automorphism group. In order to achieve similar results with
\emph{probability measure} preserving actions, we will need the
following lemma.

\begin{lemma}\label{lemma.reduction}
Suppose that we are in the setup of Lemma \ref{lemma.main}.
\begin{enumerate}
\item Suppose that $\Lambda \actson Y$ is induced from $\Lambda_1 \actson Y_1$
with $\eta(Y_1) < \infty$. Set $\Gamma_1 = \pi^{-1}(\Lambda_1)$.
Then, the orbit equivalence relation of the free, ergodic,
probability measure preserving action $\Gamma_1 \actson X \times
Y_1$ coincides with the reduction of $\cR_\pi$ to $X \times Y_1$.
\item Suppose that $\Gammatil = H * K$ where $H \cong \F_\infty \cong K$. Let
$Y_1 \subset Y$ be any measurable subset with $0 < \eta(Y_1) <
\infty$. Then, the restriction $\cR$ of $\cR_\pi$ to $X \times Y_1$
is treeable with infinite cost. Hence, by Corollary~1.2 in
\cite{hjorth}, there exists a free, ergodic, probability measure
preserving action $\F_\infty \actson X \times Y_1$ such that $\cR =
\cR(\F_\infty \actson X \times Y_1)$.
\end{enumerate}
\end{lemma}

\begin{proof}
The first part of the statement is obvious.

Since $\Gamma \cong \F_\infty$ acts freely on $X \times Y$, the
orbit equivalence relation $\cR_\pi$ is treeable.  By Proposition
II.6 in \cite{gab2}, also $\cR$ is treeable. Now, $H \cong \F_\infty
\cong K$ and $\Gamma$ is a normal subgroup of $H * K$ containing
$H$. It follows that $\Gamma$ freely splits of $H$ and that we can
choose free generators $g_n, n \in \N$ of $H$ that can be completed
into a free generating set of $\Gamma \cong \F_\infty$. Lemma II.8
in \cite{gab2} provides a treeing for $\cR$, starting from the
treeing for $\cR(\Gamma \actson X \times Y)$ given by the set of
free generators of $\Gamma$. This treeing contains for every $n \in
\N$, the automorphisms $(x,y) \mapsto (g_n \cdot x,y)$ of $X \times
Y_1$. Hence, the cost of $\cR$ is infinite.
\end{proof}

Recall from the introduction the notion of an \emph{ergodic measure} on the Borel sets of $\R$,
as well as the associated subgroup $H_\nu \subset \R$. As explained in the preliminaries,
the group $H_\nu$ can be uncountable (without being $\R$) and with prescribed Hausdorff
dimension.

\begin{theorem} \label{thm.prescribed}
Let $\nu$ be an ergodic measure on $\R$ with associated subgroup
$H_\nu \subset \R$. Let $\cG$ be any totally disconnected unimodular
locally compact group. There exists an uncountable family
$(\si_i)_{i \in I}$ of free, ergodic, probability measure preserving
actions $\si_i : \F_\infty \actson (X_i,\mu_i)$ with the following
properties.
\begin{itemize}
\item The fundamental group of $\cR(\F_\infty \actson X_i)$ equals $\exp(H_\nu)$.
\item The outer automorphism group of $\cR(\F_\infty \actson X_i)$ is isomorphic with $\cG$.
\item The actions $(\si_i)_{i \in I}$ are not stably orbit equivalent.
\item $\rL^\infty(X)$ is an HT Cartan subalgebra of $\rL^\infty(X_i) \rtimes_{\si_i}
\F_\infty$ in the sense of \cite[Definition 6.1]{P-Betti}.
\end{itemize}
\end{theorem}
\begin{proof}
We apply Lemma \ref{lemma.main} with $\Lambdatil = \Lambda$. Below,
we construct a free ergodic m.p.\ action of $\Lambda$ on an infinite
measure space $(Y,\eta)$, $\Lambda \actson (Y,\eta)$, such that
\begin{itemize}
\item the homomorphism $\module : \Centr_\Lambda(Y) \recht \Rp$
has image $\exp(H_\nu)$ and kernel isomorphic with $\cG$,
\item the action $\Lambda \actson (Y,\eta)$ is orbit equivalent to a profinite action.
\end{itemize}
Recall here that $\Centr_\Lambda(Y)$
denotes the group of non-singular automorphisms of $Y$ that commute with the
action $\Lambda \actson (Y,\eta)$. Such automorphisms $\Delta$ automatically
scale the measure $\eta$
and the scaling factor is denoted by $\module \Delta$.

We then take $\Gamma = H * K$ where $H \cong \F_\infty \cong K$.
Using Theorem \ref{thm.main-technical}, we choose a free m.p.\
$\Gamma$-action on a probability space, $\Gamma \actson (X,\mu)$,
such that $H$ acts ergodically and rigidly, and such that
$\Embed(H^\infty,\Gamma^\infty) = [\Gamma^\infty]$. Because of
Lemmas \ref{lemma.main} and \ref{lemma.reduction}, any surjective
homomorphism $\pi : \Gamma \recht \Lambda$ with $H \subset \Ker \pi$
will provide us with an action of $\F_\infty$. Since we can choose
uncountably many different $\Ker \pi$, we actually get uncountably
many non stably orbit equivalent actions of this kind.

We now construct $\Lambda \actson (Y,\eta)$ as the product of free,
ergodic, infinite measure preserving actions $\Lambda_1 \actson
(Y_1,\eta_1)$ and $\Lambda_2 \actson (Y_2,\eta_2)$, in such a way
that
\begin{itemize}
\item $\Lambda_1$ is amenable and $\Lambda_2 \actson (Y_2,\eta_2)$ is profinite,
\item the homomorphism $\module : \Centr_{\Lambda_1}(Y_1) \recht \exp(H_\nu)$
is an isomorphism of groups,
\item the elements of $\Centr_{\Lambda_2}(Y_2)$ are measure preserving and
form a group isomorphic with $\cG$.
\end{itemize}
By ergodicity of $\Lambda_i \actson (Y_i,\eta_i)$, any $\Delta \in \Centr_\Lambda(Y)$
has the form $\Delta = \Delta_1 \times \Delta_2$, so that we can safely construct
both actions separately. By amenability of $\Lambda_1$,
the action $\Lambda_1 \actson (Y_1,\eta_1)$ is orbit equivalent with a profinite action.

\subsubsection*{{\bf Part 1.} {\it Construction of the action $\Lambda_1 \actson
(Y_1,\eta_1)$.}}

Fix a countable amenable group $G$, different from $\{e\}$ and having infinite conjugacy classes.

Define for every $0 < t < 1$, the probability space
$$(Z_t,\eta_t) := \Bigl( \{0,1\} , \frac{1}{1+t} \delta_0 +
\frac{t}{1+t} \delta_1 \Bigr)^{G} \; .$$
The group $G \times G$ acts on the set $G$ by left-right translation and so,
we define the amenable group
$$G_t := \Bigl(\frac{\Z}{2 \Z} \Bigr)^{(G)} \rtimes (G \times G) \; ,$$
acting freely on $(Z_t,\eta_t)$ by non-singular automorphisms. Note
however that these automorphisms do not preserve $\eta_t$.

Let $\nu$ be an ergodic measure on $\R$ and $Q \subset \R$ a countable subgroup such that
$\nu \circ \lambda_x = \nu$ for all $x \in Q$ and such that every $Q$-invariant Borel
function on $\R$ is $\nu$-almost everywhere constant. Set $\cJ
= Q \cap (0,1)$. Define the amenable group
$$\Lambda = \bigoplus_{t \in \cJ} G_t \quad\text{and}\quad (Y_0,\eta_0)
= \prod_{t \in \cJ} (Z_t,\eta_t) \; .$$
Note that $\Lambda$ acts on $(Y_0,\eta_0)$ by non-singular automorphisms.

Define $(Y,\eta) = (Y_0 \times \R,\mu_0 \times \nu_1)$ where $\nu_1$ is given by
$d\nu_1(x) = \exp(x) d\nu(x)$. Denote by
$\om : \Lambda_1 \times Y_1 \recht Q$ the logarithm of the Radon-Nikodym $1$-cocycle.
It follows that the action
$$\Lambda \actson Y : g \cdot (y,x) = (g \cdot y, x - \om(g,y))$$
preserves the infinite measure $\eta$. The action $\Lambda \actson
Y$ is free.

Consider the subgroup $\cL$ of $\Lambda$ defined by
$$\cL = \bigoplus_{t \in \cJ} (G \times G) \; .$$
Observe that $\cL \actson (Y_0,\eta_0)$ preserves the measure $\eta_0$ and is ergodic.
Therefore, a $\Lambda$-invariant measurable function $F$ on $Y = Y_0 \times \R$,
is $\eta$-almost everywhere equal to a function only depending on the $\R$-variable. Since
the Radon-Nikodym $1$-cocycle attains all values in $Q$, the ergodicity of $\nu$
implies that $F$ is $\eta$-almost everywhere constant.

To conclude part~1, it suffices to prove that every non-singular automorphism $\Delta$
of $(Y,\eta)$ that commutes with the $\Lambda$-action, is of the form $\Delta_z :
(y,x) \mapsto (y,x+z)$ for some $z \in H_\nu$. Write
$\Delta(y,x) = (\Delta_1(y,x),\Delta_2(y,x))$.
It follows that $\Delta_2$ is $\cL$-invariant and
hence essentially independent of the $Y_0$-variable.
Hence, $\Delta_2(y,x) = \beta(x)$, $\eta$-almost everywhere,
for some non-singular automorphism $\beta$ of $(\R,\nu)$
commuting with the $Q$-action. But then the function
$\R \recht \R : z \mapsto \beta(z) - z$ is $Q$-invariant and so,
by our assumptions on $\nu$, $\nu$-almost everywhere constant.
Hence, $\beta(z) = \lambda_x(z)$ for some $x \in \R$ and $\nu$-almost every
$z \in \R$. The non-singularity of $\beta$ implies that $x \in H_\nu$.
Composing $\Delta$ with $\id \times \lambda_{-x}$, we obtain a non-singular
automorphism $\Delta$ of $(Y,\eta)$ commuting with the $\Lambda$-action and
satisfying $\Delta(y,x) \in Y_0 \times \{x\}$ almost everywhere.
So, for $\nu$-almost every $x \in \R$, we find a non-singular automorphism
$\al_x$ of $(Y_0,\eta_0)$ such that $\Delta(y,x) = (\al_x(y),x)$ almost everywhere.
But then, almost every $\al_x$ commutes with the $\Lambda$-action on $(Y_0,\eta_0)$
and in particular with the $\cL$-action. As in the proof of \cite[Theorem 5.4]{PV},
such an $\al_x$ is the identity almost everywhere. Hence, $\Delta$ is the identity
almost everywhere.

We denote as $\Lambda_1 \actson (Y_1,\eta_1)$ the action $\Lambda \actson (Y,\eta)$
constructed above.

\subsubsection*{{\bf Part 2.} {\it Construction of the action $\Lambda_2 \actson
(Y_2,\eta_2)$.}}

Choose a countable dense subgroup $\Lambda \subset \cG$ acting on $\cG$ by left multiplication.
Equip $\cG$ with its Haar measure. Let $\cK_n \subset \cG$ be a decreasing sequence of compact
open subgroups of $\cG$ with trivial intersection. Defining $Y_n = \cG / \cK_n$, the
profiniteness
of $\Lambda \actson \cG$ follows, with $\Lambda_n = \Lambda$ for all $n$.

We leave it as an exercise to prove that every non-singular automorphism $\Delta$ of $\cG$
commuting with the $\Lambda$-action is given by right multiplication with an element in $\cG$.
We denote as $\Lambda_2 \actson (Y_2,\eta_2)$ the action $\Lambda \actson (\cG,\text{Haar})$
constructed above.
\end{proof}

\begin{remark}\label{rem.extras}
The unimodularity assumption on $\cG$ in Theorem \ref{thm.prescribed} is not essential.
If $\cG$ is no longer unimodular, the image of $\cG$ under the modular function
will be part of the fundamental group and the outer automorphism group will be
given by the kernel of the modular function.

Since in Lemma \ref{lemma.main}, the action $\Lambda \actson (Y,\eta)$ is only assumed
to be orbit equivalent with a profinite action, other locally compact groups than the
totally disconnected ones can be covered in Theorem \ref{thm.prescribed}. In particular,
any second countable, locally compact, abelian group arises, since any countable dense
subgroup is still abelian and hence, amenable.
\end{remark}

The following lemma was used in the proof of Lemma \ref{lemma.main}.
Its proof is very similar to arguments in A.1 of \cite{P-Betti}. We
include the details for completeness.

\begin{lemma} \label{lem.intertwine}
Let $(A,\tau)$ be an abelian von Neumann algebra with faithful normal tracial state
$\tau$ and $0 < \eps < 1$. Let $B,C \subset A$ von Neumann subalgebras. If
$\|E_C(b) - b\|_2 \leq \eps$ for all $b$ in the unit ball of $B$, there exists
a projection $q \in A$ with $\tau(q) \geq 1-28\sqrt{\eps}$ such that $Bq \subset Cq$.
\end{lemma}
\begin{proof}
Consider the basic construction $\langle A,e_C \rangle$ of the inclusion $C \subset A$,
equipped with its natural semi-finite trace $\Tr$. Denote
$$\cK = \convex \{ u e_C u^* \mid u \in \cU(A)\} \; .$$
By our assumption $\|x - e_C \|_{2,\Tr} \leq 2\eps$ for all $x \in \cK$. Let $x$ be
the element of minimal $\| \, \cdot \, \|_{2,\Tr}$ in the weak closure of $\cK$.
Then, $\|x - e_C \|_{2,\Tr} \leq 2\eps$, $0 \leq x \leq 1$, $\Tr(x) \leq 1$ and
$x \in B' \cap \langle A,e_C \rangle$. Denote by $p$ the spectral projection $p
= \chi_{[1/2,1]}(x)$. It follows that $p$ is a projection in
$B' \cap \langle A,e_C \rangle$ such that $\|p - e_C \|_{2,\Tr} \leq 4\eps$.

Denote by $\cE_C$ the center-valued weight from $\langle A,e_C \rangle^+$ to
the extended positive part of $C$. Note that $\Tr = \tau \circ \cE_C$.

The image of the projection $p$ is a Hilbert $B$-$C$-subbimodule $K$ of $\rL^2(A,\tau)$.
As a right $C$-module, the isomorphism class of $K$ is determined by $\cE_C(p)$.
We have $\|\cE_C(p) - 1 \|_2 \leq 4 \eps$. If we view $(C,\tau)$ as $\rL^\infty(X,\mu)$,
the function $\cE_C(p)$ takes values in $\N \cup \{+\infty\}$. Defining $z$ as the
largest projection in $C$ such that $\cE_C(p) z = z$, it follows that $\|1-z\|_2 \leq
4 \eps$. Then,
$$\|p(1-z)\|_{2,\Tr}^2 = \tau(\cE_C(p) (1- z)) = \tau(\cE_C(p) - z) = \tau(\cE_C(p) - 1)
+ \tau(1 - z) \; .$$
So, $\|p(1-z)\|_{2,\Tr}^2 \leq 8 \eps$ and hence, $\|p - pz \|_{2,\Tr} \leq 3 \sqrt{\eps}$.
But then,
$$\|pz - e_C \|_{2,\Tr} \leq \|pz - p\|_{2,\Tr} + \|p - e_C \|_{2,\Tr} \leq 3 \sqrt{\eps} +
4 \eps \leq 7 \sqrt{\eps} \; .$$

Replacing $p$ by $pz$, we have found a projection $p$ in $B' \cap \langle A,e_C \rangle$
such that $\|p - e_C \|_{2,\Tr} \leq 7 \sqrt{\eps}$ and such that $\cE_C(p) = z$, where
$z$ is a projection in $C$ satisfying $\|1-z\|_2 \leq 4 \eps$.

As a $C$-module, the image $K$ of $p$ is then isomorphic to $\rL^2(Cz)$. This provides
us with $v \in \rL^2(A)z$ such that $K$ is the
$\rL^2$-closure of $vC$ and such that $E_C(v^*v) = z$. Hence, $p = v e_C v^*$. Denote
by $q$ the support projection of $v$. Because $K$ is also a $B$-module, we have
$Bq \subset Cq$. It remains to prove that $\tau(q) \geq 1-28 \sqrt{\eps}$.

We have $\|v e_C v^* -e_C \|_{2,\Tr} \leq 7 \sqrt{\eps}$, meaning that
$1 + \tau(z) - 2 \|E_C(v)\|_2^2 \leq 49 \eps$. Since $\tau(1-z) \leq 4 \eps$,
we get $1 - \|E_C(v)\|_2^2 \leq 27 \eps$. The left hand side equals $\|v - E_C(v)\|_2^2$
and so $\|v - E_C(v)\|_2 \leq 6 \sqrt{\eps}$. Also $\|v\|_2 \leq 1$, implying that
$\|v^* v - E_C(v)^* E_C(v) \|_1 \leq 12 \sqrt{\eps}$. Applying $E_C$ and using
$E_C(v^*v) = z$, we conclude that $\|v^* v - z \|_1 \leq 24 \sqrt{\eps}$.
But then, since $q \leq z$, we have $\tau(z-q) \leq 24 \sqrt{\eps}$.
We also have $\tau(1-z) \leq 4 \eps \leq 4 \sqrt{\eps}$, finally
leading to $\tau(q) \geq 1-28 \sqrt{\eps}$.
\end{proof}

\section{The II$_1$ factor case} \label{sec.vnalg}

The part of Theorem \ref{main-theorem} pertaining to the group
measure space II$_1$ factors $M_i=L^\infty(X)\rtimes_{\sigma_i} \F_\infty$ follows now
immediately from Theorem
\ref{thm.prescribed} and the results in
\cite{P-Betti} about uniqueness of HT Cartan subalgebras. Indeed, by
\cite[Theorem 6.2]{P-Betti}, if $M_1$ and $M_2$ are II$_1$ factors
with HT Cartan subalgebras $A_i \subset M_i$ and corresponding
equivalence relations $\cR_i$, any $^*$-isomorphism $\pi : M_1
\recht M_2$ can be composed with an inner automorphism Ad$u$ of
$M_2$ so that Ad$(u)(\pi(A_1))= A_2$, thus implementing an orbit
equivalence of $\cR_1$ and $\cR_2$. Similarly for stable
isomorphisms. This shows in particular that if $A \subset M$ is an
HT Cartan subalgebra with corresponding equivalence relation $\cR$,
then $\cF(M) = \cF(\cR)$ and $\Out(M) \cong \bH^1(\si) \rtimes
\Out(\cR)$.

When $\si : \F_\infty \actson (X,\mu)$ is a  free ergodic m.p.\
action on a probability space, we can say the following about
$\bH^1(\si)$ (see Lemma 3.1 in \cite{P-comp}). Denote by $\bG$ the
abelian Polish group of unitaries in $\rL^\infty([0,1])$ and by
$\bG^\infty$ the direct product of infinitely many copies of $\bG$.
Set $\bG_0 = \bG/\T$. Then, $\bZ^1(\si)$ can be identified with
$\bG^\infty \cong \bG$. Further, $\bB^1(\si) \cong \bG_0$ and the
subgroup $\bB^1(\si) \subset \bZ^1(\si)$ is closed if $\si$ is
strongly ergodic. The latter holds for all the actions $\si_i$ in
Theorem \ref{main-theorem}, because these $\si_i$ extend the natural
action $\F_2 \actson \T^2$, which is strongly ergodic by
\cite{schmidt1}.


\begin{thebibliography}{15}\setlength{\itemsep}{-1mm} \setlength{\parsep}{0mm} \small
\bibitem{aar-nad} {\sc J. Aaronson and M. Nadkarni},
$L_\infty$ eigenvalues and $L_2$ spectra of non-singular transformations. {\it Proc.
London Math. Soc.} {\bf 55} (1987), 538--570.

\bibitem{burger} {\sc M. Burger}, Kazhdan constants for $\SL(3,\Z)$,
{\it J. Reine Angew. Math.} {\bf 413} (1991), 36--67.

\bibitem{connesthesis} {\sc A. Connes}, Une classification des facteurs de
type III, {\it Ann. Ec. Norm. Sup.} {\bf 6} (1973), 133--252.

\bibitem{connes} {\sc A. Connes}, A factor of type II$_1$ with
  countable fundamental group. {\it J. Operator Theory} {\bf 4}
  (1980), 151--153.
  
\bibitem{CJ} {\sc A. Connes and V. Jones}, A II$_1$ factor with two nonconjugate Cartan subalgebras.
{\it Bull. Amer. Math. Soc. (N.S.)} {\bf 6} (1982), 211--212. 

\bibitem{fm} {\sc J. Feldman and C.C. Moore}, Ergodic equivalence
relations, cohomology, and von Neumann algebras I, II, {\it Trans.
Amer. Math. Soc.} {\bf 234} (1977), 289--324, 325--359.

\bibitem{fur1} {\sc A. Furman}, Orbit equivalence rigidity, {\it Ann.
of Math.} {\bf 150} (1999), 1083--1108.

\bibitem{fur2} {\sc A. Furman}, Outer automorphism groups of some
  ergodic equivalence relations. {\it Comment. Math. Helv.} {\bf 80} (2005), 157--196.

\bibitem{gab1} {\sc D. Gaboriau}, Invariants $l^2$ de relations
  d'\'{e}quivalence et de groupes. {\it Publ. Math. Inst. Hautes \'{E}tudes
  Sci.} {\bf 95} (2002), 93--150.

\bibitem{gab2} {\sc D. Gaboriau}, Co\^{u}t des relations d'\'{e}quivalence et
  des groupes. {\it Invent. Math.} {\bf 139} (2000), 41--98.


\bibitem{GP} {\sc D. Gaboriau and S. Popa}, An uncountable family
of nonorbit equivalent
actions of $\F_n$. {\it J. Amer. Math. Soc.} {\bf 18} (2005), 547--559.

\bibitem{ge1} {\sc S.L. Gefter}, Outer automorphism group of the ergodic
equivalence relation generated by translations of dense subgroup of compact
group on its homogeneous space. {\it Publ. Res. Inst. Math. Sci.}
{\bf 32} (1996), 517--538.


\bibitem{gg} {\sc S.L. Gefter and V.Ya. Golodets}, Fundamental
groups for ergodic actions and actions with unit fundamental groups.
{\it Publ. Res. Inst. Math. Sci.} {\bf 24} (1988), 821--847.

\bibitem{giord-skand} {\sc T. Giordano and G. Skandalis}, Krieger factors isomorphic
to their tensor square and pure point spectrum flows. {\it J. Func. Anal.} {\bf 64}
(1985), 209--226.

\bibitem{hjorth} {\sc G. Hjorth}, A lemma for cost attained.
{\it  Ann. Pure Appl. Logic}
{\bf 143}  (2006), 87--102.

\bibitem{HMP} {\sc B. Host, J.-F. M\'{e}la and F. Parreau},
Non-singular transformations and spectral analysis of measures.
{\it Bull. Soc. Math. France} {\bf 119} (1991), 33--90.


\bibitem{IPP} {\sc  A. Ioana, J. Peterson and S. Popa}, Amalgamated
  free products of $w$-rigid factors and calculation of their symmetry
  groups. {\it Acta Math.} {\bf 200} (2008), 85--153.

\bibitem{kahane-salem} {\sc J.-P. Kahane and R. Salem}, Ensembles parfaits et
s\'{e}ries trigonom\'{e}triques. Hermann, Paris, 1963.

\bibitem{kechris} {\sc A.S. Kechris},
Classical descriptive set theory. {\it Graduate Texts in Mathematics} {\bf 156}.
Springer-Verlag, New York, 1995.

\bibitem{LG} {\sc B. Le Gac}, Some properties of Borel subgroups of real numbers.
{\it Proc. Amer. Math. Soc.} {\bf 87} (1983), 677--680.

\bibitem{MNP} {\sc V. Mandrekar, M. Nadkarni and D. Patil},
Singular invariant measures on the line.
{\it Studia Math.} {\bf 35} (1970), 1--13.

\bibitem{monod-shalom} {\sc N. Monod and Y. Shalom}, Orbit equivalence
  rigidity and bounded cohomology. {\it Ann. Math.} {\bf 164} (2006), 825--878.

\bibitem{MvN1} {\sc F. Murray, J. von Neumann}, On rings of
operators, {\it Ann. of Math.} {\bf 37} (1936), 116--229.

\bibitem{MvN2} {\sc F. Murray, J. von Neumann}, Rings of operators
IV, {\it Ann. Math.} {\bf 44} (1943), 716--808.

\bibitem{NPS} {\sc R. Nicoara, S. Popa and R. Sasyk}, On II$_1$ factors
arising from $2$-cocycles of $w$-rigid groups. {\it J. Func. Anal.}
{\bf 242} (2007), 230--246.

\bibitem{OPI} {\sc N. Ozawa and S. Popa},
On a class of $\mathrm{II}_1$ factors
with at most one Cartan subalgebra. {\it Ann. Math.,} to appear. {\tt arXiv:0706.3623}

\bibitem{P-corr} {\sc S. Popa}, Correspondences, INCREST preprint No.
56/1986, \newline www.math.ucla.edu/$\sim$popa/preprints.html

\bibitem{P-Betti} {\sc S. Popa}, On a class of
type II$_1$ factors with Betti numbers invariants. {\it Ann. of
  Math.} {\bf 163} (2006), 809--899.

\bibitem{P-strong} {\sc S. Popa}, Strong rigidity of II$_1$ factors
arising from malleable actions of $w$-rigid groups II, {\it Invent.
Math.} {\bf 165} (2006), 409--452.


\bibitem{P-gap} {\sc S. Popa}, On the superrigidity of malleable
actions with spectral gap, {\it J. of Amer. Math. Soc.}, on line
Sept. 26, 2007, http://www.ams.org/journals/jams/0000-000-00/

\bibitem{P-comp} {\sc S. Popa}, Some computations of $1$-cohomology groups and
construction of non-orbit-equivalent actions. {\it J. Inst. of Math.
Jussieu} {\bf 5} (2006), 309--332.

\bibitem{PS} {\sc S. Popa and D. Shlyakhtenko},
Universal properties of $L(\F_\infty)$ in subfactor theory, {\it
Acta Mathematica}, {\bf 191} (2003), 225--257.

\bibitem{PV} {\sc S. Popa and S. Vaes}, Strong rigidity of
generalized Bernoulli actions and computations of their symmetry
groups. {\it Adv. Math.} {\bf 217} (2008), 833--872.


\bibitem{schmidt1} {\sc K. Schmidt}, Asymptotically invariant sequences
and an action of $\SL(2, \Z)$ on the $2$-sphere, {\it Israel. J.
Math.} {\bf 37} (1980), 193--208.


\bibitem{singer} {\sc I.M. Singer}, Automorphisms of finite factors.
{\it Amer. J. Math.} {\bf 177} (1955), 117--133.

\bibitem{To} {\sc A. T\"{o}rnquist}, Orbit equivalence and actions of $\F_n$.
{\it J. Symbolic Logic} {\bf 71} (2006), 265--282.

\bibitem{vaes} {\sc S. Vaes},
Explicit computations of all finite index bimodules for a family of
II$_1$ factors. {\it Ann. Sc. Ecole Norm. Sup.}, to appear. {\tt arXiv:0707.1458}



\end{thebibliography}
\end{document}